\documentclass[12pt,a4paper]{article}
\usepackage{amsmath,amssymb,comment}

\newtheorem{example}{Example}

\newcommand{\N}{{\mathbb N}}

\newcommand{\C}{{\mathbb C}}
\newcommand{\Z}{{\mathbb Z}}

\newcommand{\Dsc}{{\mathcal D}}

\newcommand{\GCD}{{\mathrm{gcd}}}
\title{Bernstein-Sato polynomials and analytic non-equivalence of plane curve singularities}
\author{Toshinori Oaku
\\Department of Mathematics, Tokyo Woman's Christian University, 
\\Suginami-ku, Tokyo, 167-8585, Japan
}
\date{October 26, 2019}
\begin{document}
\maketitle

\begin{abstract}      
We compute Bernstein-Sato polynomials of 
some pairs of topologically equivalent plane curve singularities. 
Some pairs have the same Tjurina number but  
distinct Bernstein-Sato polynomials, which implies that they are not
analytically equivalent. 
\end{abstract}

{\small
\begin{flushleft}
\end{flushleft}
}

\section{Introduction}

The purpose of this paper is to present some examples of pairs of  
plane curve singularities with the same topological type 
and the same Tjurina number that have different Bernstein-Sato polynomials, 
which implies that the pair are not analytically equivalent. 

In general, 
let $f$ and $g$ be complex analytic functions defined on a 
neighborhood of the origin $0$ in $\C^n$.  
Then $f$ and $g$, or more precisely, the germs $f=0$ and $g=0$ of 
complex hypersurfaces, are said to be analytically equivalent at the origin 
if there exist a germ of biholomorphic map $\varphi$ at $0$ 
and a germ of holomorphic function $u$ at $0$ 
such that $\varphi(0) = 0$, $u(0) \neq 0$,  and $g = u(f \circ \varphi)$. 

On the other hand, $f$ and $g$ are said to be topologically equivalent 
at the origin if there exists a homeomorphism $\varphi$ of a neighborhood 
$U$ onto $V$ such that $\varphi(0) = 0$ and 
$\varphi(\{x \in U | f(x) = 0\}) = \{x \in V | g(x) = 0\}$. 

Let $\C\{x\}$ be the ring of convergent power series in the variables 
$x = (x_1,\dots,x_n)$.  
We denote by $J_f$ the ideal of $\C\{x\}$ generated by 
$f$ and its derivatives $\partial f/\partial x_i$ ($i=1,\dots,n$). 
The quotient algebra $T(f) := \C\{x\}/J_f$ is called the Tjurina algebra 
and its dimension as the vector space over $\C$ is called 
the Tjurina number of $f$, which we denote by $\tau(f)$. 

Mather and Yau \cite{MY} proved under the condition 
$f$ and $g$ have isolated singularity at $0$ that 
$f$ and $g$ are analytically equivalent if and only if 
$T(f)$ and $T(g)$ are isomorphic as $\C$-algebras. 
In particular, the equality $\tau(f)=\tau(g)$ 
is a necessary conditon for $f$ and $g$ to be analytically equivalent. 
See \cite{KangKim}, \cite{SL} for classification of analytic equivalence of 
curves defined by some specific polynomials including cubics, 
and e.g.\ \cite{CGL} for related problems for curves. 
In any case, 
not much seems to be  known about classification of analytic
quivalence classes of germs of plane curves. 

The local Bernstein-Sato polynomial, which is also called the local 
$b$-function, 
is also invariant under analytic equivalence. 
Let us denote by $\Dsc$ the ring of linear differential operators 
with coefficients in $\C\{x\}$. 
Then the local $b$-function $b_f(s)$ of $f$ at the origin $0$ 
is the nonzero polynomial $b(s)$ of the least degree 
in an indeterminate $s$ that satisfies
\[
b(s)f^s \in \Dsc[s] f^{s+1},
\]
that is, there exists a polynomial $P(s)$ in $s$ with coefficients in $\Dsc$ 
such that $b(s)f^s = P(s)f^{s+1}$ holds. 
Here $f^s$ is regarded as a formal function, on which $\Dsc[s]$ acts 
naturally. 
The above $b_f(s)$ is uniquelly determined up to nonzero constant muliple, 
hence is unique if we impose that $b_f(s)$ be monic.

Kashiwara \cite{Kashiwara} proved that 
there always exists such $b_f(s)$ and its roots are negative rational numbers. 
Yano \cite{Yano} calculated the local Bernstein-Sato polynomials of 
a variety of examples including curves. 
An algorithm for computing $b_f(s)$ of an aribtrary polynomial $f$ was 
given by the present author \cite{OakuDuke}. See also 
\cite{NN} for improvements.

On the other hand, 
complete classification of topological equivalence is well-known 
for germs of plane curves. 
In what follows, we restrict our attention to germs of 
holomorphic functions in two variables 
which are irreducible in the unique factorization domain $\C\{x,y\}$ 
in two variables $x$, $y$. 

By a locally holomorphic change of the coordinates, the germ of 
the curve $f=0$ is parametrized by a Puiseux expansion
\begin{equation}\label{eq:par}
x = t^p,\quad 
y = t^q + c_{q+1}t^{q+1} + c_{q+2}t^{q+2} + \cdots
\end{equation}
in the complex parameter $t$ with $|t|$ sufficiently small, 
where $p,q$ are positive integers and $c_{q+1}$, 
$c_{q+2}$, $\dots$ are complex numbers.  We set $c_q = 1$. 
We may also assume $p < q$ and that the greatest common divisor 
of the elements of the set 
$\{p,q\} \cup \{j \in \Z \mid j > q,\, c_j \neq 0\}$ 
is one. 
Then the Puiseux characteristic of (\ref{eq:par}) is defined as follows 
(see \cite{Wall}):
First set
\[
r_1 = \min\{ r \in \Z \mid r \geq q,\,c_r\neq 0,\, p \not| r\}, 
\quad
e_1 = \GCD(p,r_1),
\]
where $\GCD$ means the greatest common divisor. 
If $e_1 > 1$, then set
\[
r_2 = \min\{ r \in \Z \mid r \geq q,\,c_r\neq 0,\, e_1 \not| r\}, 
\quad
e_2 = \GCD(e_1,r_2). 
\]
Define $r_i$ and $e_i$ recursively in the same way. 
Then we have $e_m=1$ for some $m$ and terminate this procedure. 
The sequence $(p; r_1,\dots,r_m)$ is called the Puiseux characteristic 
of (\ref{eq:par}). 
In paricular, if $p$ and $q$ are relatively prime, then the Puiseux characteristic is simply $(p;q)$. 

It is a classical result attributed to Burau and Zariski 
dating back to the 1930s that 
the Puiseux characteristic is in one-to-one correspondence to each 
topological equivalence class of the germs of plane curves  
(see e.g., \cite{Wall}).  
In particular, if $p$ and $q$ are relatively prime, then the curve germ  
defined by (\ref{eq:par}) is topologically equivalent to that of 
$y^p - x^q = 0$. 
It should be noted that the Alexander polynomial of the knot defined by 
(\ref{eq:par}), which can be thought of as a topological counterpart 
of the Bernstein-Sato polynomial, 
plays an essential role in the proof of the classical theorem above. 

Yano \cite{Yano2} 
made a conjecture about the generic Bernstein-Sato polynomial of a plane curve germ in terms of its Puiseux characteristic. 
There are many works related to his conjecture;  
see e.g., \cite{CN},\cite{ABC},\cite{Blanco}. 
However, complete theoretical description of the behavior 
of the Bernstein-Sato polynomials of the  curves 
with the same Puiseux characteristic seems to be unknown.

\section{Examples}

We give examples of topologically equivalent germs of plane curves 
some of which have the same Tjurina number but differernt 
Bernstein-Sato polynomials. 

In what follows $p$ and $q$ are relatively prime positive integers 
with $p < q$ and  $f_0(x,y) = y^p -x^q$.  
Let $\N$ be the set of non-negative integers and 
set 
\[
G(p,q) = \{i \in \N \mid i > q\} \setminus (\N p + \N q).
\] 
Then it is easy to see that by a holomorphic change of local coordinates, 
the parametrization (\ref{eq:par}) can be transformed to a simple form
\[
x = t^p,\quad 
y = t^q + \sum_{r\in G(p,q)} c_{r}t^{r}. 
\]
Among such parametrizations, we pick up the one 
\begin{equation}\label{eq:par2}
x = t^p,\quad 
y = t^q + t^r
\end{equation}
for each $r \in G(p,q)$. Note that if $c_r \neq 0$, then this is 
analytically equivalent to 
$x = t^p$, $y = t^q + c_rt^r$. 

For each $r \in G(p,q)$, let $f_r(x,y)$ be the polynomial whose 
germ at $0$ is the defining function of the plane curve germ parametrized by 
(\ref{eq:par2}). 
See 2.3 of \cite{Wall} for a method of computing $f_r$, other than 
the elimination method based on an appropriate Gr\"obner basis. 
We denote by $\tau(f_r)$ the Tjurina number of $f_r$ and 
by $b_r(s)$ the local Bernstein polynomial of $f_r$ at $0$. 
Note that $f_0$ and $f_r$ with $r\in G(p,q)$ are all topologically equivalent. 
Note also that there is an explicit formula 
(see 6.4 of \cite{KashiwaraBook}) 
for the Bernstein-Sato polynomials 
of quasi-homogeneous polynomials with isolated singularity, which applies to 
the binomial $x^p - y^q$.

The following examples were computed by using the library file 
``nn\_ndbf.rr'' of Risa/Asir developed by Nishiyama and Noro \cite{NN}.

\begin{example}\rm 
We set $p=4$, $q=9$. Then we have $f_0 = y^4-x^9$ and
 $G(4,9) = \{10,11,14,15,19,23 \}$. 
Corresponding polynomials are
\begin{align*}
f_{10} &= y^4 - 2x^5y^2 - 4x^7y - x^9 + x^{10},
\\
f_{11} &= y^4 - 4x^5y^2 - x^9 + 2x^{10} - x^{11},
\\
f_{14} &= y^4 - 2x^7y^2 - 4x^8y - x^9 + x^{14},
\\
f_{15} &=  y^4 - 4x^6y^2 - x^9 + 2x^{12} - x^{15},
\\
f_{19} &= y^4 - 4x^7y^2 - x^9 + 2x^{14} - x^{19},
\\
f_{23} &= y^4 - 4x^8y^2 - x^9 + 2x^{16} - x^{23}.
\end{align*}
The Tjurina numbers of $f_r$ are as follows:
\[
\begin{array}{|c|c|c|c|c|c|c|c|c|}
\hline
f & f_0 & f_{10} & f_{11} & f_{14} & f_{15} & f_{19} & f_{23} 
\\ \hline 
\tau(f) & 24 & 21 & 21 & 23 & 22 & 23 & 24 
\\ \hline
\end{array}
\]
The local Bernstein-Sato polynomial $b_0(s) = b_{f_0}(s)$ of $f_0$ at $0$ is 
\begin{align*}
b_0(s) &= 
(36s+13)(36s+17)(36s+25)(36s+29)(36s+31)(36s+35)
\\&\times
(36s+37)(36s+41)(36s+43)(36s+47)(36s+55)(36s+59)
\\&\times
(18s+11)(18s+13)(18s+17)(18s+19)(18s+23)(18s+25)
\\&\times
(12s+7)(12s+11)(12s+13)(12s+17)(6s+5)(6s+7)(s+1). 
\end{align*}
%
%
The Bernstein-Sato polynomials of $f_r$ with $r \in G(4,9)$ are as 
follows: 
\begin{align*}
b_{10}(s) &= b_0(s)\frac{(36s+19)(36s+23)(18s+7)(12s+5)}
                        {(36s+55)(36s+59)(18s+25)(12s+17)},
\\
b_{11}(s) &= b_0(s)\frac{(36s+19)(36s+23)(12s+5)}
                        {(36s+55)(36s+59)(12s+17)},
\\
b_{14}(s) &= b_0(s)\frac{36s+23}
                        {36s+59},
\\
b_{15}(s) &= b_0(s)\frac{(36s+19)(36s+23)}
                        {(36s+55)(36s+59)},
\\
b_{19}(s) &= b_{14}(s),
\\
b_{23}(s) &= b_{0}(s).
\end{align*}
Thus $f_{10}$ and $f_{11}$ have the same Tjurina number $21$ but have
different Bernstein-Sato polynomials. 
We do not know if each of the pairs $(f_0,f_{23})$ and 
$(f_{14},f_{19})$ is one with analytically equivalent germs. 
\end{example}

\begin{example}\rm 
$p=5$, $q=6$, $f_0 = y^5-x^6$, $G(5,6) = \{7,8,9,13,14,19\}$. 
The corresponding polynomials are
\begin{align*}
f_{7} &= y^5 - 5x^4y^2 - 5x^5y - x^6 - x^7,
\\
f_{8} &= y^5 - 5x^4y^2 - 5x^6y - x^6 - x^8,
\\
f_{9} &= y^5 - 5x^3y^3 + 5x^6y - x^6 - x^9,
\\
f_{13} &= y^5 - 5x^5y^2 - 5x^9y - x^6 - x^{13},
\\
f_{14} &= y^5 - 5x^4y^3 + 5x^8y - x^6 - x^{14},
\\
f_{19} &= y^5 - 5x^5y^3 + 5x^{10} - x^6 - x^{19}.
\end{align*}
The Tjurina numbers are as follows:
\[
\begin{array}{|c|c|c|c|c|c|c|c|}
\hline
f & f_0 & f_7 & f_8 & f_9 & f_{13} & f_{14} & f_{19}
\\ \hline 
\tau(f) & 20 & 18 & 18 & 18 & 20 & 19 & 20 
\\ \hline
\end{array}
\]
The Bernstein-Sato polynomial of $f_0$ is
\begin{align*}
b_0(s) &= 
(30s+11)(30s+17)(30s+23)(30s+29)(30s+31)(30s+37)
\\&\times
(30s+43)(30s+49)(15s+8)(15s+11)(15s+13)(15s+14)
\\&\times
(15s+16)(15s+17)(15s+19)(15s+22)(10s+7)(10s+9)
\\&\times
(10s+11)(10s+13)(s+1). 
\end{align*}
Those of $f_r$ are as follows:
\begin{align*}
b_{7}(s) &= b_0(s)\frac{(30s+13)(30s+19)(15s+7)}{(30s+43)(30s+49)(15s+22)},
\\
b_{8}(s) &= b_0(s)\frac{(30s+13)(30s+19)}{(30s+43)(30s+49)},
\\
b_{9}(s) &= b_0(s)\frac{(30s+19)(15s+7)}{(30s+49)(15s+22)},
\\
b_{13}(s) &= b_0(s),
\\
b_{14}(s) &= b_0(s)\frac{30s+19}{30s+49},
\\
b_{19}(s) &= b_0(s). 
\end{align*}
Thus $f_7$, $f_8$, $f_9$ have the same Tjurina number $18$ but 
distinct Bernstein-Sato polynomials. 
\end{example}

\begin{example}\rm 
$p=5$, $q=7$, $f_0 = y^5-x^7$, $G(5,7) = \{8,9,11,13,16,18,23 \}$. 
The corresponding polynomials are
\begin{align*}
f_8 &= y^5 - 5x^3y^3 + 5x^6y - x^7 - x^8,
\\
f_9 &= y^5 - 5x^5y^2 - 5x^6y - x^7 - x^9,
\\
f_{11} &= y^5 - 5x^5y^2 - 5x^8y - x^7 - x^{11},
\\
f_{13} &= y^5 - 5x^4y^3 + 5x^8y - x^7 - x^{13},
\\
f_{16} &= y^5 - 5x^6y^2 - 5x^{11} - x^7 - x^{16}, 
\\
f_{18} &= y^5 - 5x^5y^3 + 5x^{10}y - x^7 - x^{18},
\\
f_{23} &= y^5 - 5x^6y^3 + 5x^{12}y - x^7 - x^{23}.
\end{align*}
The Tjurina numbers are as follows: 
\[
\begin{array}{|c|c|c|c|c|c|c|c|c|c|}
\hline
f & f_0 & f_8 & f_9 & f_{11} & f_{13} & f_{16} & f_{18} & f_{23} 
\\ \hline 
\tau(f) & 24 & 21 & 22 & 22 & 22 & 24 & 23 & 24  
\\ \hline
\end{array}
\]
The Bernstein-Sato polynomial of $f_0$ is
\begin{align*}
b_0(s) &= 
(35s+12)(35s+17)(35s+19)(35s+22)(35s+24)(35s+26)
\\&\times
(35s+27)(35s+29)(35s+31)(35s+32)(35s+33)(35s+34)
\\&\times
(35s+36)(35s+37)(35s+38)(35s+39)(35s+41)(35s+43)
\\&\times
(35s+44)(35s+46)(35s+48)(35s+51)(35s+53)(35s+58)(s+1)
.
\end{align*}
Those of other $f_r$ are as follows:
\begin{align*}
b_8(s) &= \frac{(35s+13)(35s+18)(35s+23)}{(35s+48)(35s+53)(35s+58)},
\\
b_9(s) &= \frac{(35s+16)(35s+18)(35s+23)}{(35s+51)(35s+53)(35s+58)},
\\
b_{11}(s) &= \frac{(35s+16)(35s+23)}{(35s+51)(35s+58)},
\\
b_{13}(s) &= \frac{(35s+18)(35s+23)}{(35s+53)(35s+58)},
\\
b_{16}(s) &= b_0(s),
\\
b_{18}(s) &= \frac{35s+23}{35s+58},
\\
b_{23}(s) &= b_0(s).
\end{align*}
Thus $f_9$, $f_{11}$, and $f_{13}$ have the same Tjurina number $22$ 
but distinct Bernstein-Sato polynomials. 
\end{example}


\noindent
\textbf{Acknowledgement}

This work was supported in part by JSPS Grant-in-Aid for Scientific Research (C) 26400123. 

\end{document}